\documentclass[11pt]{article}
\usepackage[letterpaper, margin=1in]{geometry}
\usepackage[utf8]{inputenc}
\usepackage[T1]{fontenc}
\usepackage{amsmath,amssymb,amsfonts,amsthm,mathtools}
\usepackage{graphicx}
\usepackage{booktabs}
\usepackage{siunitx}
\usepackage{authblk}
\usepackage[labelfont=bf,format=plain,justification=raggedright,singlelinecheck=false]{caption}
\usepackage[colorlinks=true,linkcolor=blue,citecolor=blue,urlcolor=blue]{hyperref}
\usepackage[section]{placeins} % keep floats within their section (no drift into References)

\DeclareMathOperator{\eig}{eig}

\usepackage{aliascnt}
\newtheorem{thm}{Theorem}[section]
\newaliascnt{cor}{thm}\newtheorem{cor}[cor]{Corollary}\aliascntresetthe{cor}
\newaliascnt{lm}{thm}\newtheorem{lm}[lm]{Lemma}\aliascntresetthe{lm}
\newaliascnt{df}{thm}\newtheorem{df}[df]{Definition}\aliascntresetthe{df}
\newaliascnt{remrk}{thm}\aliascntresetthe{remrk}
\newaliascnt{prp}{thm}\newtheorem{prp}[prp]{Proposition}\aliascntresetthe{prp}
\usepackage[capitalise]{cleveref}

\crefname{thm}{Theorem}{Theorems}
\Crefname{thm}{Theorem}{Theorems}
\crefname{prp}{Proposition}{Propositions}
\Crefname{prp}{Proposition}{Propositions}
\crefname{lm}{Lemma}{Lemmas}
\Crefname{lm}{Lemma}{Lemmas}
\crefname{cor}{Corollary}{Corollaries}
\Crefname{cor}{Corollary}{Corollaries}
\crefname{df}{Definition}{Definitions}
\Crefname{df}{Definition}{Definitions}

\title{Probabilistic Analysis of Least Squares, Orthogonal Projection, and QR Factorization Algorithms Subject to Gaussian Noise}
\author[1]{Ali Lotfi}
\author[2]{Julien Langou}
\author[3]{Mohammad Meysami}
\affil[1]{Department of Computer Science, University of Saskatchewan, Saskatoon, SK S7N 5A2, Canada}
\affil[2]{Department of Applied Mathematics, University of Colorado Denver, Denver, Colorado 80201, USA}
\affil[3]{Department of Mathematics, The University of Tulsa, Tulsa, OK 74104, USA}

\begin{document}
\maketitle

\begin{abstract}
{\color{black}
We consider the effect of Gaussian perturbations on least-squares residuals, orthogonal projections, and QR-type algorithms. The problem that motivated our investigations is as follows: suppose that a full column-rank matrix \(B\in\mathbb{R}^{m\times n}\) has already been computed, and suppose that a new normalized column \(q=(x+y)/\|x+y\|_2\) is to be appended to \(B\), where \(x\perp\operatorname{span}(B)\) is the ideal orthogonal component and \(y\) represents the orthogonalization error. How large can the condition number \(\kappa([B,q])\) of the resulting matrix \([B,q]\) become? While we provide a Weyl-type bound on the singular values of \([B,q]\), in terms of the extremal singular values of \(B\) and the quantity \(\|B^T y\|_2/\|x+y\|_2\), we also derive exact probability laws for norms and projection residuals under Gaussian perturbations. Finally, we use these probability laws to derive probabilistic condition-number bounds for QR-type processes with imperfect orthogonalization and exact normalization.
}
\end{abstract}

\section{Introduction}
\subsection{History}
{\color{black}

Stability analyses have been conducted for many of the QR algorithms that have been proposed. For algorithms based on the Modified Gram-Schmidt algorithm, for example, the loss of orthogonality has been studied in \cite{bjorck1992loss,paige2006modified,pereira2022numerical,lotfi2022numerical}. Related Gram-Schmidt based algorithms, such as the GSXO algorithm, have been studied in \cite{lotfi2022numerical}. These studies typically consider the stability of the algorithm in the context of machine precision and worst-case rounding-error models. The analyses often make assumptions regarding the number of columns of the matrix, the conditioning of the input matrix, and the unit roundoff. However, in practice, algorithms like MGS can often produce useful orthonormal bases even outside of the most conservative regimes covered by these estimates. Thus, a complementary probabilistic analysis of the loss of orthogonality in algorithms like MGS is motivated.

A second class of literature considers the condition number and singular values of random or randomly perturbed matrices \cite{chen2005condition,ratnarajah2004eigenvalues,edelman1988eigenvalues,rudelson2008littlewood}. Related approaches based on smoothed analysis and perturbations are discussed in \cite{tao2010smooth,vu2007condition,szarek1991condition}. Abstract approaches to incorporating probabilistic methods into considerations of numerical analysis are discussed in \cite{burgisser2008probability}. Bounding methods for the singular values of random matrices, such as those discussed in \cite{sengupta1999distributions}, lead to analyses of the condition number of those matrices. Each of these methods is relevant to the topic of this work, but does not directly answer the question investigated here regarding how imperfect orthogonalization affects the conditioning of the computed QR factor when a new column is appended.
}

\subsection{Overview}
The present paper studies a tractable Gaussian model for the
orthogonality-damaging part of the error. The model is not intended to be a
complete floating-point error model. It is a first-order model for the
accumulated effect of many small coefficient errors in an orthogonalization
step. In a Gram-Schmidt step, computed coefficients have the form
\(\widehat r_{jk}=r_{jk}+\delta_{jk}\), and the computed residual contains a
perturbation of the form
\[
-\sum_{j<k}\delta_{jk}q_j.
\]
Each \(\delta_{jk}\) is produced by a sum of many elementary arithmetic errors.
Under the usual independence or weak-dependence heuristic for accumulated small
errors, a Gaussian approximation is natural. The results below should therefore
be read as a probabilistic perturbation analysis for orthogonalization error,
not as a replacement for deterministic backward-error analysis.

The analysis is organized around three elementary objects: the norm of a vector
subject to Gaussian noise, an orthogonal projection subject to Gaussian noise,
and a normalized least-squares residual subject to Gaussian noise. These objects
then feed into condition-number bounds for QR-type algorithms in which
normalization is accurate but orthogonalization is imperfect.

\begin{figure}[htbp]
\centering
\includegraphics[width=\textwidth]{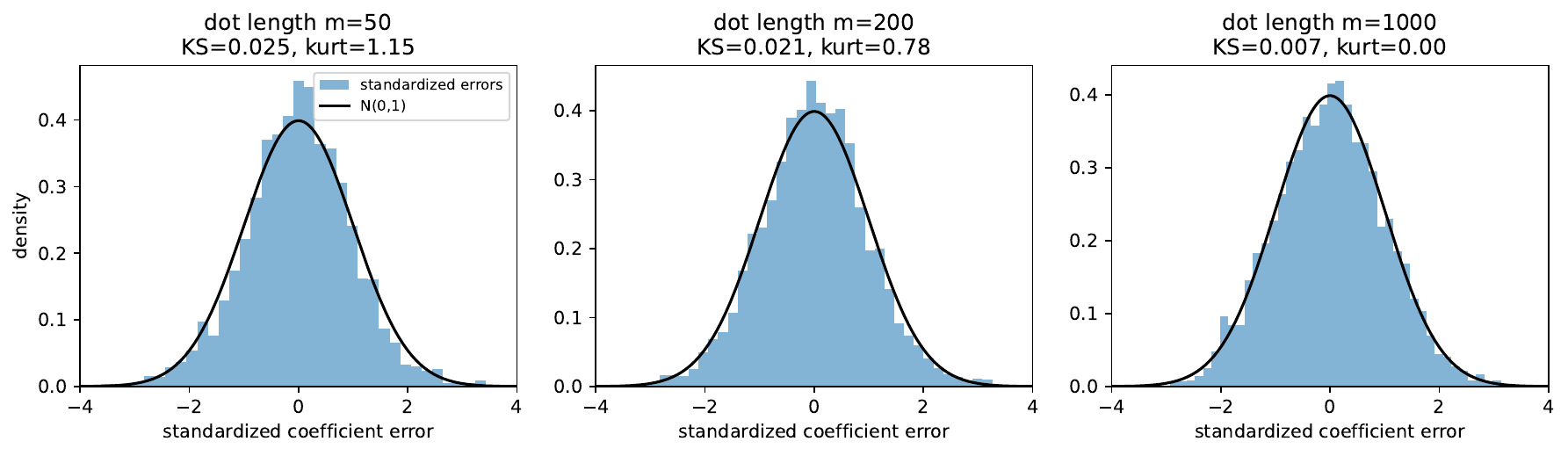}
\caption{Standardized single-precision inner-product errors compared with
\(N(0,1)\). As the vector length \(m\) increases, the empirical distribution
approaches a standard normal (decreasing Kolmogorov--Smirnov distance and excess
kurtosis), consistent with a central-limit explanation for the accumulation of
elementary rounding errors in a Gram--Schmidt coefficient.}
\label{fig:clt-normality}
\end{figure}

These modeling assumptions can be checked numerically, and they are consistent
with recent probabilistic rounding-error analyses that treat rounding errors as
mean-zero, mean-independent random variables \cite{higham2019new,connolly2021stochastic},
including QR-specific studies of Householder \cite{connolly2023probabilistic}
and modified Gram--Schmidt \cite{zou2024probabilistic}. First, the central-limit
viewpoint is supported directly: standardized single-precision inner-product
errors approach a standard normal as the vector length \(m\) grows
(\Cref{fig:clt-normality}), with the Kolmogorov--Smirnov distance to
\(N(0,1)\) falling from \(0.025\) at \(m=50\) to \(0.007\) at \(m=1000\) and the
excess kurtosis from \(1.15\) to essentially \(0\). Second, the perturbation
scale can be tied to the unit roundoff: a one-parameter span-confined model with
per-column standard deviation \(\sigma_k = c\,u\,\|a_k\|_2\), with a single
constant \(c\approx0.4\) calibrated once, reproduces the median loss of
orthogonality \(\kappa(\widehat Q)-1\) of single-precision modified
Gram--Schmidt to within a small factor across several configurations; it
underpredicts the upper tail, so it should be read as a typical-scale rather
than a worst-case model. A companion damage score
\(u\,\|a_k\|_2\sqrt{k-1}/r_{kk}\) correlates with the observed column overlap
\(\|\widehat Q_{1:k-1}^T\widehat q_k\|_2\) at Spearman \(0.90\) over more than
\(14{,}000\) columns, identifying later columns and small residual norms as the
vulnerable ones. These diagnostics are reproducible from the bundled script and
are intended as motivation, not as deterministic rounding-error theorems.

\subsection{Motivation and problem-solving strategy}
This work builds on the relation between condition numbers and least-squares
residuals established by Liesen, Rozlo\v{z}n\'{i}k, and Strako\v{s}
\cite{liesen2002least}. Their result shows that the condition number of the
augmented matrix \([Q,\gamma c]\), where \(Q\) has orthonormal columns, is
determined by the least-squares residual of \(c\) from \(\operatorname{span}(Q)\).
Our contribution is to combine this deterministic viewpoint with Gaussian noise
models for residuals and projections.

More precisely, let \(B\in\mathbb{R}^{m\times n}\) have full column rank and
suppose a new normalized column has the form
\[
q=\frac{x+y}{\|x+y\|_2},\qquad x\perp \operatorname{span}(B).
\]
Here \(x\) is the ideal orthogonal component and \(y\) is the
orthogonalization error. The deterministic part of the paper bounds
\(\kappa([B,q])\) in terms of the singular values of \(B\) and the scalar
quantity \(\|B^Ty\|_2/\|x+y\|_2\). The probabilistic part evaluates the
likelihood of the corresponding projection and least-squares quantities being
large or small under Gaussian noise.

\subsection{Paper organization}
\Cref{seq:notation} fixes notation. \Cref{condition-grow-general} proves a
deterministic append-one-column condition-number bound. \Cref{sec:norm-of-noised-vector}
gives the exact distribution of \(\|X+Y\|_2\) under Gaussian noise. The
special-function representation uses the modified Bessel and Marcum-Q functions;
standard references include \cite{gil2014algorithm,baricz2009tight,kapinas2009monotonicity,simon2008digital}.
\Cref{sec:projection} applies the norm result to orthogonal projections.
\Cref{sec:ls} derives a noncentral-F law for a normalized least-squares residual.
\Cref{sec:qr-noise} combines these ingredients to obtain QR-type
condition-number bounds.

\section{Notation}\label{seq:notation}
For a matrix \(A\), \(\sigma_{\max}(A)\) and \(\sigma_{\min}(A)\) denote the largest and smallest singular values of \(A\). Thus
\[
\sigma_{\min}(A)=\sqrt{\lambda_{\min}(A^TA)}
\]
when \(A\) has full column rank. The spectral condition number is
\[
\kappa(A)=\frac{\sigma_{\max}(A)}{\sigma_{\min}(A)}.
\]
The spectral and Frobenius norms are denoted by \(\|A\|_2\) and \(\|A\|_F\), respectively. For a full column-rank matrix \(B\), \(B^\dagger\) denotes the Moore--Penrose pseudoinverse,
\[
B^\dagger=(B^TB)^{-1}B^T,
\]
so that \(\|B^\dagger\|_2=1/\sigma_{\min}(B)\).

For a subspace \(\mathcal S\subseteq\mathbb{R}^m\), \(P_{\mathcal S}\) denotes the orthogonal projector onto \(\mathcal S\), and \(P_{\mathcal S^\perp}=I-P_{\mathcal S}\). If \(B\in\mathbb{R}^{m\times n}\) has full column rank, then
\[
P_B=B(B^TB)^{-1}B^T,\qquad P_B^\perp=I-P_B.
\]
If \(Q\in\mathbb{R}^{m\times n}\) has orthonormal columns, then
\[
P_Q=QQ^T,\qquad P_Q^\perp=I-QQ^T.
\]
The notation \(N(0,\sigma^2I_m)\) denotes the centered Gaussian distribution in \(\mathbb{R}^m\) with covariance \(\sigma^2I_m\).
In deterministic append-column statements we use lower-case \(x,y\).
In probabilistic statements we use \(X\) for the fixed signal vector and
\(Y\) for the Gaussian perturbation. In sequential QR statements,
\(x_i\) denotes the predictable ideal residual at step \(i\), while
\(y_i\) denotes the random perturbation at that step.

\section{Condition Number Growth Induced by Adding a Column}\label{condition-grow-general}
We first prove a rank-two eigenvalue lemma used to control the off-diagonal perturbation in the Gram matrix.

\begin{lm}[Rank-two bordered matrix]\label{lm:rank2}
Let \(a\in\mathbb{R}^{n-1}\), \(b\in\mathbb{R}\), and
\[
A=\begin{bmatrix}0&a\\ a^T&b\end{bmatrix}.
\]
Then
\[
\eig(A)\subseteq\left\{0,\frac{b\pm\sqrt{b^2+4\|a\|_2^2}}{2}\right\}.
\]
\end{lm}

\begin{proof}
If \(a=0\), then \(A=\operatorname{diag}(0,\ldots,0,b)\), so the conclusion is immediate. Assume \(a\neq0\). Define
\[
u=\begin{bmatrix}a/\|a\|_2\\0\end{bmatrix},\qquad v=e_n.
\]
The subspace \(\operatorname{span}\{u,v\}\) is invariant under \(A\), and every vector orthogonal to both \(u\) and \(v\) lies in the nullspace of \(A\). Therefore all nonzero eigenvalues of \(A\) are eigenvalues of the restriction to \(\operatorname{span}\{u,v\}\). In the orthonormal basis \(\{u,v\}\), this restriction is
\[
\begin{bmatrix}0&\|a\|_2\\ \|a\|_2&b\end{bmatrix}.
\]
The characteristic polynomial of this \(2\times2\) matrix is
\[
\det\begin{bmatrix}-\lambda&\|a\|_2\\ \|a\|_2&b-\lambda\end{bmatrix}
=\lambda^2-b\lambda-\|a\|_2^2.
\]
Thus the nonzero eigenvalues are
\[
\lambda=\frac{b\pm\sqrt{b^2+4\|a\|_2^2}}{2}.
\]
\end{proof}

The lemma controls the off-diagonal coupling in the Gram matrix of \([B,q]\).
  The next result converts that control into a bound on \(\kappa([B,q])\), in
  terms of the singular values of \(B\) and the single scalar \(\gamma\|B^Ty\|_2\).
\begin{thm}[Append-column condition-number bound]\label{thm:cool}
Let \(B\in\mathbb{R}^{m\times n}\) have full column rank, with \(m>n\). Let \(x,y\in\mathbb{R}^m\), assume
\[
x\in\operatorname{span}(B)^\perp,\qquad x+y\neq0,
\]
and define
\[
q=\frac{x+y}{\|x+y\|_2},\qquad \gamma=\frac{1}{\|x+y\|_2}.
\]
Then
\[
\kappa([B,q])\le
\sqrt{\frac{\max\{\sigma_{\max}^2(B),1\}+\gamma\|B^Ty\|_2}{\min\{\sigma_{\min}^2(B),1\}-\gamma\|B^Ty\|_2}},
\]
provided
\[
\gamma\|B^Ty\|_2<\min\{\sigma_{\min}^2(B),1\}.
\]
If, in addition, the columns of \(B\) have unit norm and
\[
\gamma\|B^Ty\|_2\leq \delta<\sigma_{\min}^2(B),
\]
then
\[
\kappa([B,q])\le
\kappa(B)\sqrt{1+\frac{\delta(1+\kappa(B)^{-2})}{\sigma_{\min}^2(B)-\delta}}.
\]
\end{thm}

\begin{proof}
Since \(x\in\operatorname{span}(B)^\perp\), we have \(B^Tx=0\). Hence
\[
B^Tq=\frac{B^T(x+y)}{\|x+y\|_2}=\gamma B^Ty.
\]
The Gram matrix of \([B,q]\) is
\[
[B,q]^T[B,q]=
\begin{bmatrix}B^TB&\gamma B^Ty\\ \gamma y^TB&1\end{bmatrix}=H_1+H_2,
\]
where
\[
H_1=\begin{bmatrix}B^TB&0\\0&1\end{bmatrix},\qquad
H_2=\begin{bmatrix}0&\gamma B^Ty\\ \gamma y^TB&0\end{bmatrix}.
\]
The eigenvalues of \(H_1\) are the eigenvalues of \(B^TB\) together with \(1\). Therefore
\[
\lambda_{\max}(H_1)=\max\{\sigma_{\max}^2(B),1\},\qquad
\lambda_{\min}(H_1)=\min\{\sigma_{\min}^2(B),1\}.
\]
By \Cref{lm:rank2}, the eigenvalues of \(H_2\) are contained in \(\{0,\pm\gamma\|B^Ty\|_2\}\). Hence
\[
\lambda_{\max}(H_2)=\gamma\|B^Ty\|_2,\qquad
\lambda_{\min}(H_2)=-\gamma\|B^Ty\|_2.
\]
Weyl's inequalities imply
\[
\sigma_{\max}^2([B,q])\le \max\{\sigma_{\max}^2(B),1\}+\gamma\|B^Ty\|_2
\]
and
\[
\sigma_{\min}^2([B,q])\ge \min\{\sigma_{\min}^2(B),1\}-\gamma\|B^Ty\|_2.
\]
The denominator is positive by assumption, so division gives the first bound.

Now assume the columns of \(B\) have unit norm. Then for any coordinate vector \(e_j\), \(\|Be_j\|_2=1\), so
\[
\sigma_{\min}(B)\le1\le\sigma_{\max}(B).
\]
Thus the first bound reduces to
\[
\kappa([B,q])\le \sqrt{\frac{\sigma_{\max}^2(B)+\delta}{\sigma_{\min}^2(B)-\delta}}.
\]
Finally,
\[
\frac{\sigma_{\max}^2(B)+\delta}{\sigma_{\min}^2(B)-\delta}
=\kappa(B)^2\left(1+\frac{\delta(1+\kappa(B)^{-2})}{\sigma_{\min}^2(B)-\delta}\right),
\]
which proves the second bound.
\end{proof}

The bound is governed by the single quantity \(\gamma\|B^Ty\|_2\), the coupling
between the appended column and \(\operatorname{span}(B)\). When this coupling
vanishes, the Gram matrix becomes block diagonal:
\[
[B,q]^T[B,q]=
\begin{bmatrix}
B^TB&0\\
0&1
\end{bmatrix}.
\]
Thus the exact condition number is
\[
\frac{\max\{\sigma_{\max}(B),1\}}{\min\{\sigma_{\min}(B),1\}}.
\]
In the unit-column QR setting, where
\(\sigma_{\min}(B)\le 1\le \sigma_{\max}(B)\), this reduces to \(\kappa(B)\).
The probabilistic sections below estimate how large the coupling and the
corresponding residual quantities become under Gaussian noise.
\subsection{Append-column sensitivity to noise}
\begin{figure}[htbp]
\centering
\includegraphics[width=0.7\textwidth]{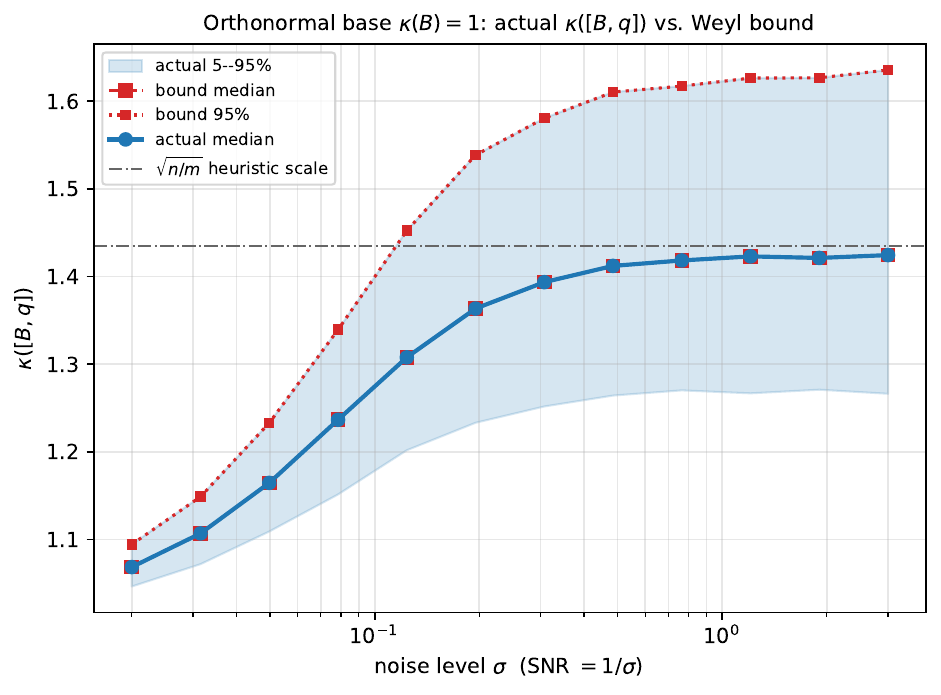}
\caption{Condition number \(\kappa([B,q])\) under increasing Gaussian perturbation (\Cref{thm:cool}), for an orthonormal base \(\kappa(B)=1\). The actual condition number rises with \(\sigma\) toward the random-vector regime, whose heuristic scale \(\sqrt{(1+\sqrt{n/m})/(1-\sqrt{n/m})}\) is dash-dotted; the Weyl-type bound coincides with the observed condition number to machine precision. The band is \(5\)--\(95\%\) over \(3{,}000\) trials per noise level.}
\label{fig:append-kappa-vs-noise}
\end{figure}

\Cref{fig:append-kappa-vs-noise} shows how \(\kappa([B,q])\) behaves as the perturbation \(y\) grows, for a base with orthonormal columns (\(\kappa(B)=1\)); here \(x\in\operatorname{span}(B)^\perp\) with \(\|x\|_2=1\), so the signal-to-noise ratio is \(1/\sigma\), and each point aggregates \(3{,}000\) trials. For this orthonormal-base experiment, the bound of \Cref{thm:cool} coincides with the observed condition number to machine precision. As the noise dominates the signal, \(q\) approaches a random unit direction; then \(\|Q^Tq\|_2\) is typically of size \(\sqrt{n/m}\), giving the heuristic condition-number scale \(\sqrt{(1+\sqrt{n/m})/(1-\sqrt{n/m})}\) (this is a typical scale, not a pointwise limit for each sample). No trial satisfying the theorem's positivity hypothesis violated the bound. For a base that is already ill-conditioned, the condition number may be dominated by the pre-existing conditioning of \(B\). When the perturbation is small, the new column is nearly orthogonal to \(\operatorname{span}(B)\) and contributes a singular value near \(1\), which lies between \(\sigma_{\min}(B)<1<\sigma_{\max}(B)\). In that regime \(\kappa([B,q])\) changes little as \(\sigma\) varies. At larger noise levels, the positivity condition \(\gamma\|B^Ty\|_2<\min\{\sigma_{\min}^2(B),1\}\) may fail, and the theorem should only be evaluated on trials for which this hypothesis holds.

This experiment supports the deterministic part of the argument: the condition number is controlled by the coupling of the new column with the existing span. In the next sections, we will quantify the probability of favorable norm and residual
events under Gaussian noise.

\section{Probabilistic Analysis of the Norm of a Vector Subject to Gaussian Noise}\label{sec:norm-of-noised-vector}

This section gives the exact distribution of the Euclidean norm of a
  Gaussian-perturbed vector. The law is expressed through the modified Bessel and
  Marcum-Q functions, recalled first.
\begin{df}[Modified Bessel function]
For a real number \(\nu\), the modified Bessel function of the first kind of order \(\nu\) is
\[
\mathcal I_\nu(t)=\sum_{j\ge0}\frac{(t/2)^{2j+\nu}}{j!\,\Gamma(\nu+j+1)}.
\]
\end{df}
The Bessel function enters the density of the noncentral chi distribution; the
  Marcum-Q function below is the corresponding survival function.
\begin{df}[Generalized Marcum-Q function]
For \(\beta\ge0\) and \(\alpha,M>0\), define
\[
\mathcal Q_M^{ma}(\alpha,\beta)=\frac{1}{\alpha^{M-1}}\int_\beta^\infty x^M\exp\left(-\frac{x^2+\alpha^2}{2}\right)\mathcal I_{M-1}(\alpha x)\,dx.
\]
For \(\alpha=0\), we use the standard limiting definition
\[
\mathcal Q_M^{ma}(0,\beta)=\frac{\Gamma(M,\beta^2/2)}{\Gamma(M)}.
\]
\end{df}
With these definitions, the norm of a Gaussian-perturbed vector has the following
closed-form survival probability.
\begin{prp}[Gaussian norm law]\label{prp:Malcum-Bessel}
Let \(X\in\mathbb{R}^m\), and let \(Y\sim N(0,\sigma^2I_m)\). Then
\[
\mathbb P(\|X+Y\|_2>\varepsilon)=\mathcal Q_{m/2}^{ma}\left(\frac{\|X\|_2}{\sigma},\frac{\varepsilon}{\sigma}\right).
\]
\end{prp}

\begin{proof}
Set \(Z=X+Y\). Then \(Z_i/\sigma\sim N(X_i/\sigma,1)\), and
\[
\frac{\|X+Y\|_2^2}{\sigma^2}=\sum_{i=1}^m\left(\frac{Z_i}{\sigma}\right)^2.
\]
This sum has a noncentral chi-square distribution with \(m\) degrees of freedom and noncentrality parameter
\[
\lambda=\sum_{i=1}^m\left(\frac{X_i}{\sigma}\right)^2=\frac{\|X\|_2^2}{\sigma^2}.
\]
The survival function of this noncentral chi-square variable is the generalized Marcum-Q function with order \(m/2\), giving the formula.
\end{proof}

The centered case is the pure-noise version of the same law.

\begin{cor}\label{cor:golden2}
If \(Y\sim N(0,\sigma^2I_m)\), then
\[
\mathbb P(\|Y\|_2\le\varepsilon)=1-\mathcal Q_{m/2}^{ma}\left(0,\frac{\varepsilon}{\sigma}\right).
\]
\end{cor}

\begin{proof}
Set \(X=0\) in \Cref{prp:Malcum-Bessel} and take the complement of the event \(\|Y\|_2>\varepsilon\).
\end{proof}

\Cref{cor:golden2} is useful when the perturbation itself, rather than a
perturbed deterministic vector, is the object of interest.

\subsection{Monte Carlo check for the Gaussian norm law}
\begin{figure}[htbp]
\centering
\includegraphics[width=0.95\textwidth]{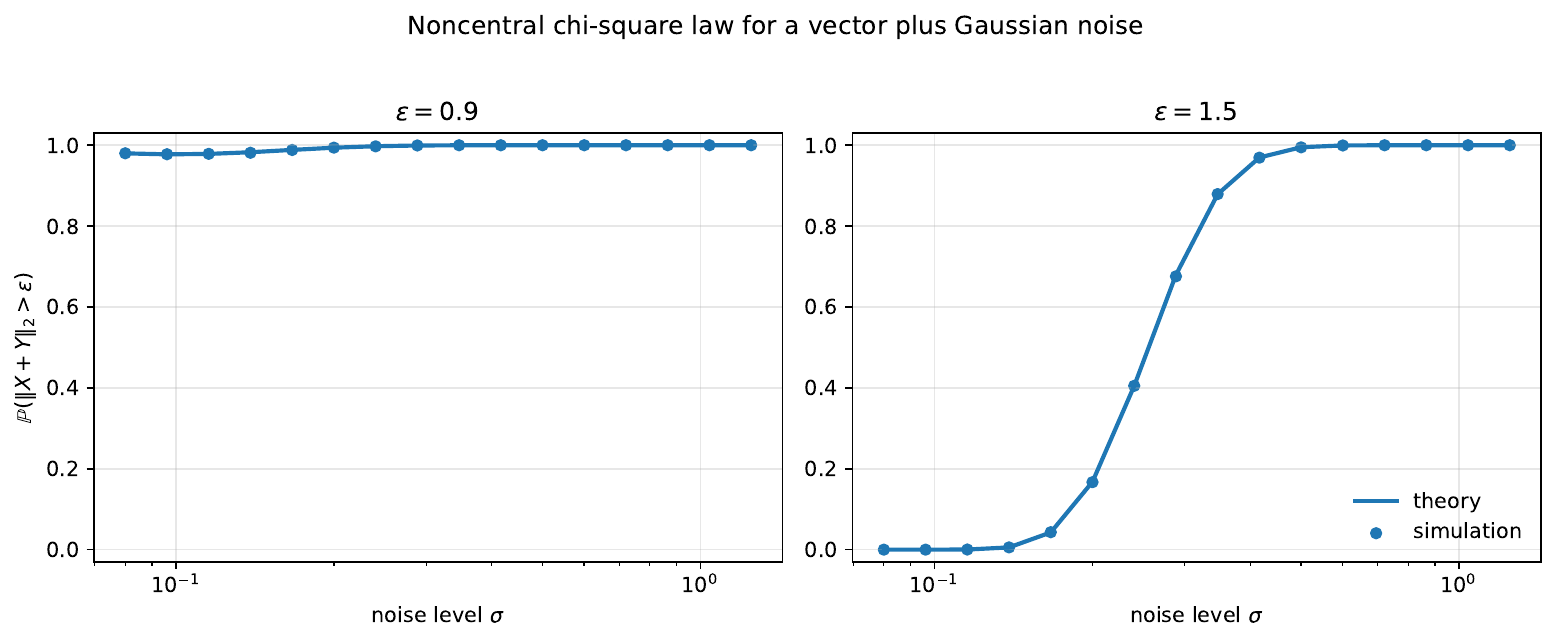}
\caption{Validation of \Cref{prp:Malcum-Bessel}. Markers are Monte Carlo estimates; curves are exact noncentral chi-square probabilities.}
\label{fig:norm-noise-validation}
\end{figure}

\Cref{fig:norm-noise-validation} compares Monte Carlo probabilities with the exact noncentral chi-square survival function for \(m=20\), \(\|X\|_2=1\), and two thresholds.

The agreement in \Cref{fig:norm-noise-validation} confirms the basic
noncentral chi-square law used throughout the rest of the paper. We next apply
the same law after projecting onto an orthogonal complement.
\section{Orthogonal Projection Subject to Gaussian Noise}\label{sec:projection}

Projecting a Gaussian-perturbed vector onto the orthogonal complement of
  \(\operatorname{span}(Q)\) lowers the ambient dimension from \(m\) to \(m-n\)
  but leaves the noise isotropic. The norm law of the previous section therefore
  applies verbatim in dimension \(m-n\).
  
\begin{prp}[Projected Gaussian norm law]\label{prp:projected-noise}
Let \(Q\in\mathbb{R}^{m\times n}\) have orthonormal columns, \(m>n\), and define
\[
P_Q=QQ^T,\qquad P_Q^\perp=I-QQ^T.
\]
Let \(N\in\mathbb{R}^{m\times(m-n)}\) have orthonormal columns spanning \(\operatorname{range}(P_Q^\perp)\), so that \(NN^T=P_Q^\perp\). Let \(X\in\operatorname{range}(P_Q^\perp)\), and let \(Y\sim N(0,\sigma^2I_m)\). Then
\[
N^TY\sim N(0,\sigma^2I_{m-n}),
\]
and
\[
\mathbb P(\|X+P_Q^\perp Y\|_2>\varepsilon)=\mathcal Q_{(m-n)/2}^{ma}\left(\frac{\|X\|_2}{\sigma},\frac{\varepsilon}{\sigma}\right).
\]
\end{prp}

\begin{proof}
Since \(N\) has orthonormal columns and \(Y\sim N(0,\sigma^2I_m)\), the vector \(N^TY\) is Gaussian with mean zero and covariance
\[
\operatorname{Cov}(N^TY)=N^T(\sigma^2I_m)N=\sigma^2I_{m-n}.
\]
Thus \(N^TY\sim N(0,\sigma^2I_{m-n})\). Since \(X\in\operatorname{range}(N)\), there exists \(\bar X\in\mathbb{R}^{m-n}\) such that \(X=N\bar X\), and \(\|\bar X\|_2=\|X\|_2\). Also
\[
P_Q^\perp Y=NN^TY=N(N^TY).
\]
Therefore
\[
\|X+P_Q^\perp Y\|_2=\|N(\bar X+N^TY)\|_2=\|\bar X+N^TY\|_2.
\]
The result follows from \Cref{prp:Malcum-Bessel} in dimension \(m-n\).
\end{proof}
The same conclusion can be stated for an arbitrary full column-rank matrix by
using its orthogonal projector.

\begin{cor}\label{cor:golden3}
Let \(B\in\mathbb{R}^{m\times n}\) have full column rank, with \(m>n\), and define
\[
P_B^\perp=I-B(B^TB)^{-1}B^T.
\]
If \(Y\sim N(0,\sigma^2I_m)\), then
\[
\mathbb P(\|P_B^\perp Y\|_2\le\varepsilon)=1-\mathcal Q_{(m-n)/2}^{ma}\left(0,\frac{\varepsilon}{\sigma}\right).
\]
\end{cor}

\begin{proof}
Take an orthonormal basis \(Q\) for \(\operatorname{span}(B)\). Then \(P_B^\perp=P_Q^\perp\) and \(X=0\) in \Cref{prp:projected-noise}.
\end{proof}

\Cref{cor:golden3} is the form used later when the previously computed columns
are represented by a full-rank matrix rather than by an explicitly orthonormal
basis.

\subsection{Monte Carlo check for projected noise}
\begin{figure}[htbp]
\centering
\includegraphics[width=0.72\textwidth]{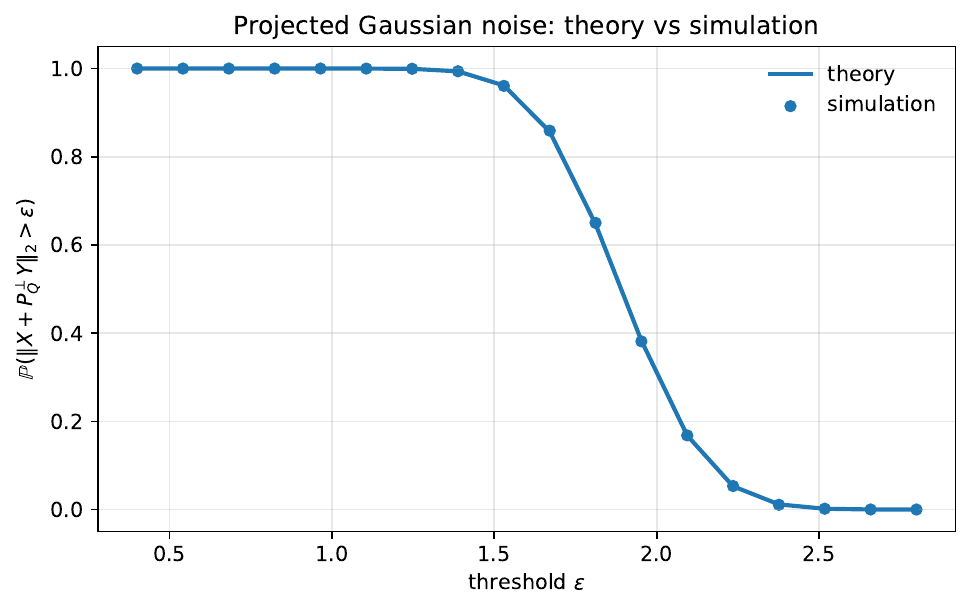}
\caption{Validation of \Cref{prp:projected-noise}. Markers are Monte Carlo estimates and the curve is the exact noncentral chi-square survival probability.}
\label{fig:projection-validation}
\end{figure}

\Cref{fig:projection-validation} validates the projected-noise formula by simulating directly in the \((m-n)\)-dimensional orthogonal complement.
Thus projection preserves the Gaussian structure while replacing the ambient
dimension \(m\) by the complementary dimension \(m-n\). The next section uses
this decomposition to study a normalized residual ratio.

\section{Least Squares Subject to Gaussian Noise}\label{sec:ls}
The next theorem describes a normalized residual ratio. This ratio is the key probability input in the later QR condition-number bounds.

\begin{thm}[Normalized residual ratio under Gaussian noise]\label{thm:noise-ls}
Let \(Q\in\mathbb{R}^{m\times n}\) have orthonormal columns, \(m>n\), and define
\[
P_Q=QQ^T,\qquad P_Q^\perp=I-QQ^T.
\]
Let \(X\in\operatorname{range}(P_Q^\perp)\), and let \(Y\sim N(0,\sigma^2I_m)\). Define
\[
\widehat r=\frac{P_Q^\perp(X+Y)}
{\sqrt{\|P_Q^\perp(X+Y)\|_2^2+\|P_QY\|_2^2}}.
\]
Then, for every \(0<\tau<1\),
\[
\mathbb P(\|\widehat r\|_2\ge\tau)=1-F'_{m-n,n}\left(\frac{n\tau^2}{(m-n)(1-\tau^2)};\frac{\|X\|_2^2}{\sigma^2}\right),
\]
where \(F'_{m-n,n}(\cdot;\lambda)\) is the CDF of the noncentral \(F\)-distribution with numerator degrees of freedom \(m-n\), denominator degrees of freedom \(n\), and numerator noncentrality parameter \(\lambda\).
\end{thm}

\begin{proof}
Set
\[
U=\frac{\|P_Q^\perp(X+Y)\|_2^2}{\sigma^2},\qquad
V=\frac{\|P_QY\|_2^2}{\sigma^2}.
\]
The two Gaussian components \(P_Q^\perp Y\) and \(P_QY\) are independent because they are orthogonal projections of an isotropic Gaussian vector. Therefore
\[
U\sim\chi'^2_{m-n}\left(\frac{\|X\|_2^2}{\sigma^2}\right),\qquad V\sim\chi^2_n,
\]
and \(U\) and \(V\) are independent. Moreover,
\[
\|\widehat r\|_2^2=\frac{U}{U+V}.
\]
Thus
\[
\|\widehat r\|_2\ge\tau
\Longleftrightarrow
\frac{U}{U+V}\ge\tau^2
\Longleftrightarrow
(1-\tau^2)U\ge\tau^2V.
\]
Equivalently,
\[
\frac{U/(m-n)}{V/n}\ge\frac{n\tau^2}{(m-n)(1-\tau^2)}.
\]
The ratio \(U/(m-n)\) divided by \(V/n\) has the stated noncentral \(F\)-distribution. Therefore the desired probability is the upper tail of this distribution.
\end{proof}

\subsection{Monte Carlo check for the residual ratio}
\begin{figure}[htbp]
\centering
\includegraphics[width=0.72\textwidth]{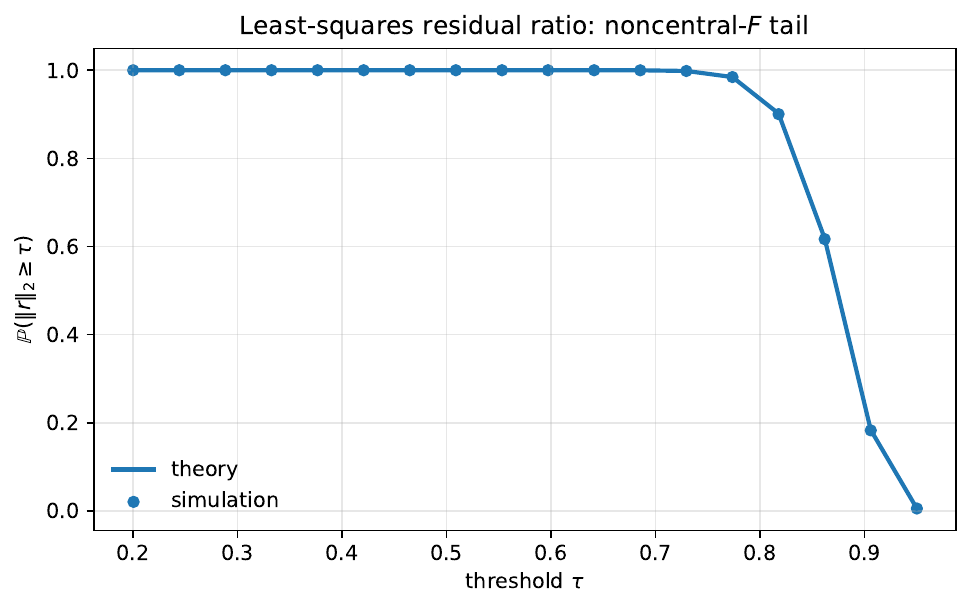}
\caption{Validation of the least-squares noise theorem. The theory is the noncentral-\(F\) upper tail; markers are Monte Carlo estimates.}
\label{fig:least-squares-validation}
\end{figure}

\Cref{fig:least-squares-validation} compares the empirical probability \(\mathbb P(\|\widehat r\|_2\ge\tau)\) with the noncentral-\(F\) upper tail.
The noncentral-\(F\) law is the probability input for the QR condition-number
bounds. The next section combines it with deterministic estimates for appending
a column.

\section{QR Factorization with Imperfect Orthogonalization}\label{sec:qr-noise}
We consider QR-type processes in which each computed column is normalized, but
the preceding orthogonalization step is perturbed.

\begin{prp}[Separating previous conditioning from a new column]\label{prp:min-max-cond}
Assume \(B\in\mathbb{R}^{m\times n}\) has full column rank, let \(c\in\mathbb{R}^m\), let \(\gamma>0\), and let \(B=QR\) be the thin QR factorization of \(B\), where \(Q^TQ=I_n\). Then:
\begin{enumerate}
\item\label{prp:min-max-cond-1}
\[
\sigma_{\max}([B,\gamma c])\le \max\{\sigma_{\max}(B),\gamma\}\,\sigma_{\max}([Q,c]).
\]
\item\label{prp:min-max-cond-2}
\[
\sigma_{\min}([B,\gamma c])\ge \min\{\sigma_{\min}(B),\gamma\}\,\sigma_{\min}([Q,c]).
\]
\item\label{prp:min-max-cond-3}
If \([Q,c]\) has full column rank, then
\[
\kappa([B,\gamma c])\le \frac{\max\{\sigma_{\max}(B),\gamma\}}{\min\{\sigma_{\min}(B),\gamma\}}\,\kappa([Q,c]).
\]
Equivalently,
\[
\kappa([B,\gamma c])\le \max\left\{\kappa(B),\frac{\|B\|_2}{\gamma},\gamma\|B^\dagger\|_2\right\}\kappa([Q,c]).
\]
Also, if \([Q,\gamma c]\) has full column rank, then
\[
\kappa([B,\gamma c])\le \frac{\max\{\sigma_{\max}(B),1\}}{\min\{\sigma_{\min}(B),1\}}\,\kappa([Q,\gamma c]).
\]
Equivalently,
\[
\kappa([B,\gamma c])\le \max\left\{\kappa(B),\|B\|_2,\|B^\dagger\|_2\right\}\kappa([Q,\gamma c]).
\]
\end{enumerate}
\end{prp}

\begin{proof}
Since \(B=QR\),
\[
[B,\gamma c]=[Q,c]\begin{bmatrix}R&0\\0&\gamma\end{bmatrix}.
\]
The singular values of \(R\) are the singular values of \(B\). Hence the block diagonal factor has largest singular value \(\max\{\sigma_{\max}(B),\gamma\}\) and smallest singular value \(\min\{\sigma_{\min}(B),\gamma\}\). The standard inequalities
\[
\sigma_{\max}(AD)\le\sigma_{\max}(A)\sigma_{\max}(D),\qquad
\sigma_{\min}(AD)\ge\sigma_{\min}(A)\sigma_{\min}(D)
\]
prove items \ref{prp:min-max-cond-1} and \ref{prp:min-max-cond-2}. Dividing proves the first condition-number bound.

Let \(S=\sigma_{\max}(B)\) and \(s=\sigma_{\min}(B)\). Then
\[
\frac{\max\{S,\gamma\}}{\min\{s,\gamma\}}=\max\left\{\frac{S}{s},\frac{S}{\gamma},\frac{\gamma}{s}\right\},
\]
which gives the displayed equivalent form. The alternative bound follows from
\[
[B,\gamma c]=[Q,\gamma c]\begin{bmatrix}R&0\\0&1\end{bmatrix}
\]
and the same singular-value argument.
\end{proof}

The factor \(\kappa([Q,c])\) isolates the effect of the appended column on an
  already-orthonormal base. The next result evaluates that factor through the
  least-squares residual of the new column.
  
\begin{prp}[Residual-to-condition transfer]\label{prp:cond-ls}
Assume \(B\in\mathbb{R}^{m\times n}\) has full column rank and columns of norm \(1\). Let \(B=QR\) be its thin QR factorization. Let \(c\in\mathbb{R}^m\), let \(\gamma>0\), and set
\[
r=\min_z\|c-Bz\|_2=\min_z\|c-Qz\|_2,\qquad \alpha=1+\gamma^2\|c\|_2^2.
\]
If \(r>0\), then
\[
\kappa([B,\gamma c])\le \kappa(B)\frac{\alpha+\sqrt{\alpha^2-4\gamma^2r^2}}{2\gamma r}.
\]
Equivalently, if \(\rho=\min_z\|\gamma c-Bz\|_2=\gamma r\), then
\[
\kappa([B,\gamma c])\le \kappa(B)\frac{\alpha+\sqrt{\alpha^2-4\rho^2}}{2\rho}.
\]
\end{prp}

\begin{proof}
Because the columns of \(B\) have unit norm, \(\sigma_{\min}(B)\le1\le\sigma_{\max}(B)\). Hence
\[
\frac{\max\{\sigma_{\max}(B),1\}}{\min\{\sigma_{\min}(B),1\}}=\kappa(B).
\]

By the alternative bound in item~\ref{prp:min-max-cond-3} of
\Cref{prp:min-max-cond},
\[
\kappa([B,\gamma c])\le\kappa(B)\kappa([Q,\gamma c]).
\]

The residual norm of \(c\) from \(\operatorname{span}(B)\) is the same as its residual norm from \(\operatorname{span}(Q)\). Applying Theorem~2.3 of \cite{liesen2002least} to the orthonormal matrix \(Q\) gives
\[
\kappa([Q,\gamma c])=\frac{\alpha+\sqrt{\alpha^2-4\gamma^2r^2}}{2\gamma r}.
\]
Multiplication by \(\kappa(B)\) proves the first formula. The second formula is the same identity with \(\rho=\gamma r\).
\end{proof}

For a normalized appended column, the preceding proposition simplifies to a
one-parameter residual bound.

\begin{cor}[Unit-column residual bound]\label{cor:cond-ls}
Assume \(B\in\mathbb{R}^{m\times n}\) has full column rank and columns of norm \(1\). Let \(q\in\mathbb{R}^m\) satisfy \(\|q\|_2=1\), and define the scalar residual norm
\[
\rho=\min_z\|q-Bz\|_2.
\]
If \(\rho>0\), then
\[
\kappa([B,q])\le \kappa(B)\frac{1+\sqrt{1-\rho^2}}{\rho}.
\]
\end{cor}

\begin{proof}
Apply \Cref{prp:cond-ls} with \(c=q\) and \(\gamma=1\). Since
\(\|q\|_2=1\), we have
\[
\alpha=1+\gamma^2\|q\|_2^2=2.
\]
The scalar residual \(r\) in \Cref{prp:cond-ls} is exactly
\[
r=\min_z\|q-Bz\|_2=\rho.
\]
Therefore
\[
\kappa([B,q])
\le
\kappa(B)
\frac{2+\sqrt{4-4\rho^2}}{2\rho}
=
\kappa(B)\frac{1+\sqrt{1-\rho^2}}{\rho}.
\]
\end{proof}

Thus \(\rho\) measures how much of the appended unit vector remains outside
\(\operatorname{span}(B)\); larger \(\rho\) gives a smaller condition-number
multiplier.

\Cref{prp:cond-ls} and \Cref{cor:cond-ls} show that the condition number is controlled
once the residual norm \(\rho\) is known. Combined with the residual law of
\Cref{thm:noise-ls}, this gives the probabilistic statement below.
  
\begin{cor}[Single-step probabilistic QR bound]\label{cor-main}
Let \(B\in\mathbb{R}^{m\times n}\) have full column rank and unit-norm columns, with \(m>n\), and let \(B=QR\) be its thin QR factorization. Let \(X\in\operatorname{range}(P_Q^\perp)\), let \(Y\sim N(0,\sigma^2I_m)\), and define
\[
q=\frac{X+Y}{\|X+Y\|_2}.
\]
Then, for every \(0<\tau<1\), if \(\|P_Q^\perp q\|_2\ge\tau\), then
\[
\kappa([B,q])\le \kappa(B)\frac{1+\sqrt{1-\tau^2}}{\tau}.
\]
Moreover,
\[
\mathbb P(\|P_Q^\perp q\|_2\ge\tau)=1-F'_{m-n,n}\left(\frac{n\tau^2}{(m-n)(1-\tau^2)};\frac{\|X\|_2^2}{\sigma^2}\right).
\]
\end{cor}

\begin{proof}
Since \(q\) has unit norm and \(\operatorname{span}(B)=\operatorname{span}(Q)\), the scalar residual norm of \(q\) from \(\operatorname{span}(B)\) is
\[
\rho=\|P_Q^\perp q\|_2.
\]
If \(\rho\ge\tau\), then \((1+\sqrt{1-\rho^2})/\rho\le(1+\sqrt{1-\tau^2})/\tau\) for \(0<\tau\le\rho\le1\). Applying \Cref{cor:cond-ls} proves the condition-number bound. The probability formula is exactly \Cref{thm:noise-ls} applied to the normalized residual ratio.
\end{proof}

Iterating the one-step result gives an unconditional multi-step bound when the
ideal residual norms are bounded below deterministically.

\begin{cor}[Multi-step product bound]\label{cor:qr-noise-product}
Assume a QR-type process produces normalized columns
\[
\widehat q_i=\frac{x_i+y_i}{\|x_i+y_i\|_2},\qquad x_i\in\operatorname{span}(\widehat Q_{1:i-1})^\perp,
\]
and assume that \(\widehat Q_{1:i-1}\) has full column rank for every \(i=2,\ldots,n\).
Suppose that \(x_i\) is \(\mathcal F_{i-1}\)-measurable and that, conditional on the previous columns (the filtration \(\mathcal F_{i-1}\) generated by \(y_2,\ldots,y_{i-1}\)), the perturbation \(y_i\sim N(0,\sigma^2I_m)\) is independent of the past. Let \(\tau_i\in(0,1)\) for \(i=2,\ldots,n\), and assume there are deterministic constants \(\ell_i>0\) with
\[
\|x_i\|_2\ge\ell_i,\qquad i=2,\ldots,n.
\]
Then, with probability at least
\[
1-\sum_{i=2}^n F'_{m-i+1,i-1}\left(\frac{(i-1)\tau_i^2}{(m-i+1)(1-\tau_i^2)};\frac{\ell_i^2}{\sigma^2}\right),
\]
we have
\[
\kappa(\widehat Q)\le\prod_{i=2}^n\frac{1+\sqrt{1-\tau_i^2}}{\tau_i}.
\]
\end{cor}

\begin{proof}
For each \(i\ge2\), apply \Cref{cor-main} with \(B=\widehat Q_{1:i-1}\); the columns of this matrix have unit norm by construction. Let
\[
E_i=\{\|P_{\widehat Q_{1:i-1}}^\perp\widehat q_i\|_2\ge\tau_i\}.
\]
On \(E_i\),
\[
\kappa(\widehat Q_{1:i})\le \kappa(\widehat Q_{1:i-1})\frac{1+\sqrt{1-\tau_i^2}}{\tau_i},
\]
so multiplying over \(i=2,\ldots,n\) gives the product bound on \(\bigcap_iE_i\). Conditional on \(\mathcal F_{i-1}\), \Cref{cor-main} gives
\[
\mathbb P(E_i^c\mid\mathcal F_{i-1})=F'_{m-i+1,i-1}\left(\frac{(i-1)\tau_i^2}{(m-i+1)(1-\tau_i^2)};\frac{\|x_i\|_2^2}{\sigma^2}\right).
\]
The noncentral-\(F\) cumulative distribution function is nonincreasing in its noncentrality parameter. This monotonicity follows from the stochastic monotonicity of the noncentral chi-square distribution in its noncentrality parameter, since the numerator of the noncentral \(F\) ratio shifts to the right as the noncentrality increases. Therefore the lower bound \(\|x_i\|_2\ge\ell_i\) yields
\[
\mathbb P(E_i^c\mid\mathcal F_{i-1})\le F'_{m-i+1,i-1}\left(\frac{(i-1)\tau_i^2}{(m-i+1)(1-\tau_i^2)};\frac{\ell_i^2}{\sigma^2}\right),
\]
a deterministic bound. Taking expectations and applying the union bound to \(\bigcup_iE_i^c\) gives the stated probability.
\end{proof}

\subsection{Probability--conditioning tradeoffs}

The corollaries above combine the deterministic residual-to-condition estimate
from \Cref{cor:cond-ls} with the noncentral-\(F\) residual law of
\Cref{thm:noise-ls}. The independent append-column bound of
\Cref{thm:cool} is validated separately in
\Cref{fig:append-kappa-vs-noise}.

\begin{figure}[htbp]
\centering
\includegraphics[width=0.66\textwidth]{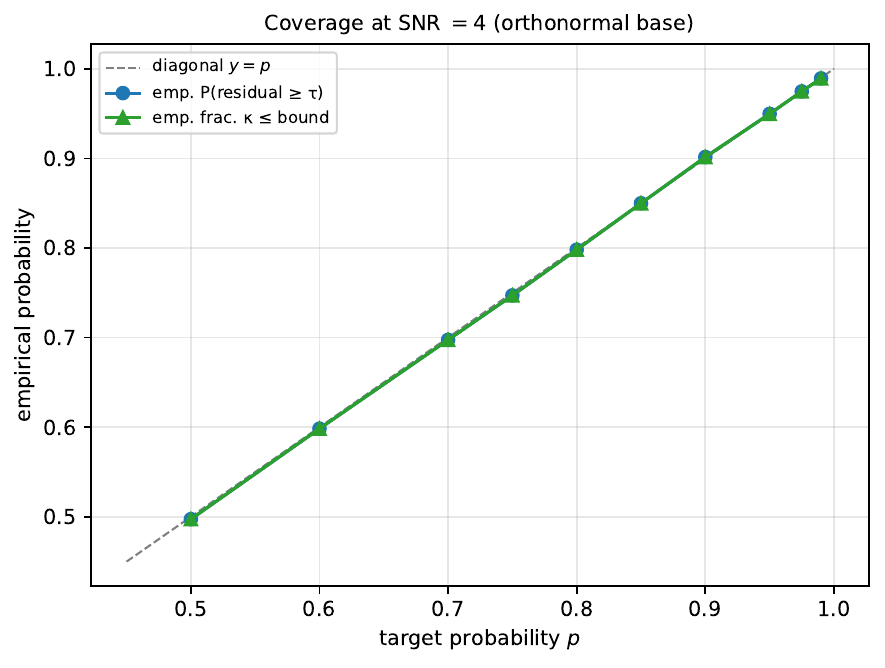}
\caption{Coverage of \Cref{cor-main} at signal-to-noise ratio
\(\|X\|_2/\sigma=4\) on an orthonormal base. In this setting, the
noncentral-\(F\) law predicts the residual exceedance probability, and this
matches the empirical fraction of sampled matrices whose condition number falls
below the certified bound.}
\label{fig:append-coverage}
\end{figure}

For the single-step probabilistic bound, we choose \(\tau\) from the
noncentral-\(F\) law for a target probability \(p\), and then compare the
certified condition-number multiplier with Monte Carlo samples of
\(\kappa([B,q])\). \Cref{fig:append-coverage} shows one representative
calibration slice at signal-to-noise ratio \(\|X\|_2/\sigma=4\).

\begin{figure}[htbp]
\centering
\includegraphics[width=\textwidth]{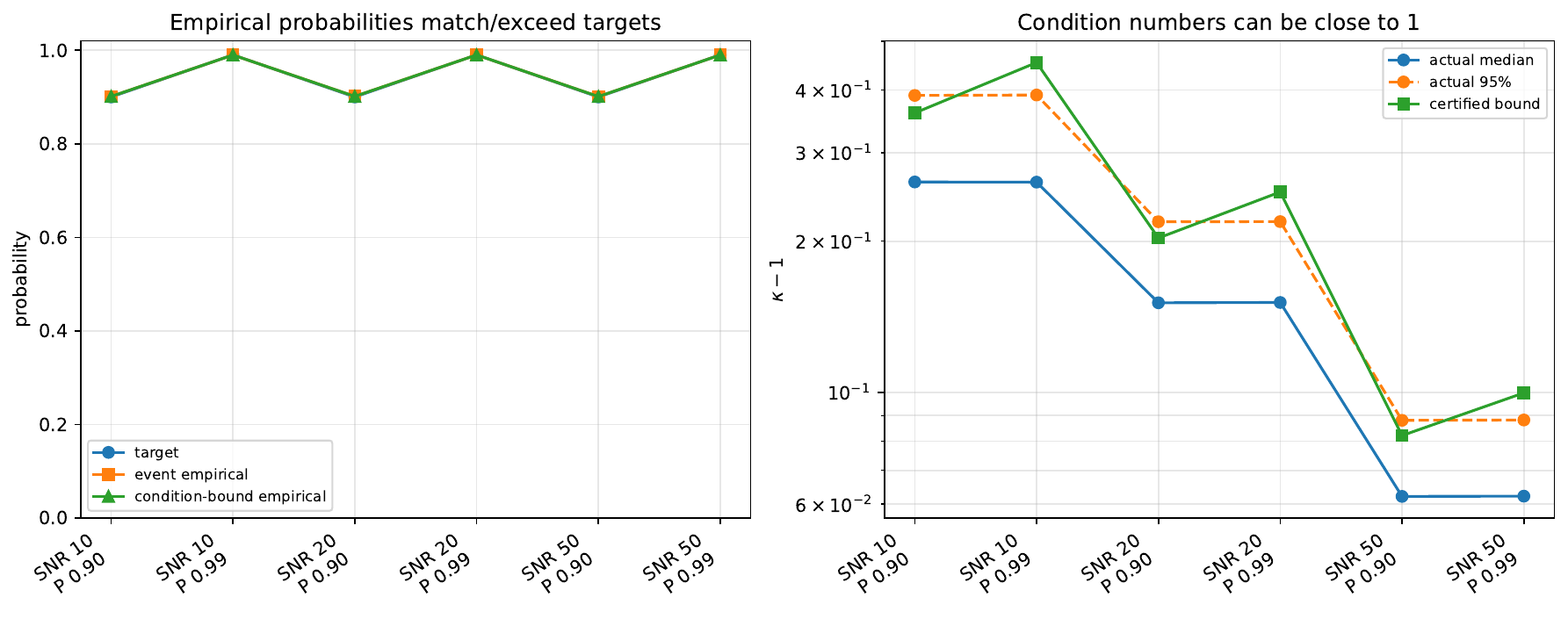}
\caption{Empirical validation of the single-append bound
(\Cref{cor-main}). Left: the empirical exceedance probability matches
the target probability. Right: the certified condition-number multiplier captures
the Monte Carlo condition numbers, with the sharpest behavior in the
orthonormal-base case.}
\label{fig:single-step-empirical}
\end{figure}

\begin{figure}[htbp]
\centering
\includegraphics[width=0.66\textwidth]{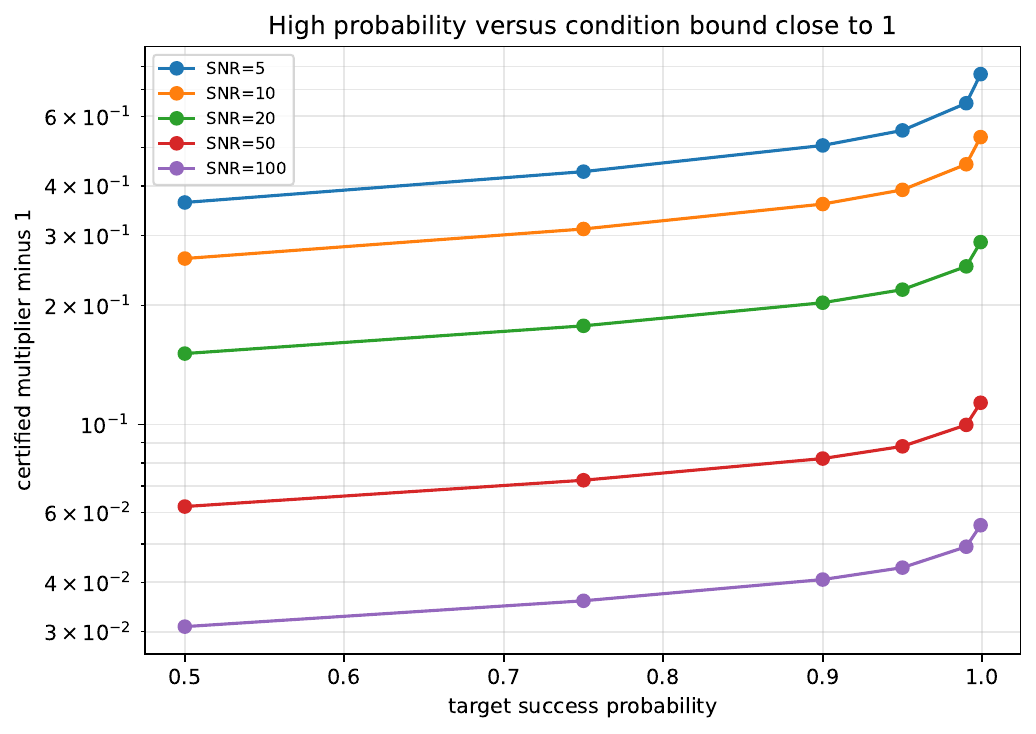}
\caption{Single-append tradeoff between a target success probability and the
certified condition-number multiplier of \Cref{cor-main}. A higher
signal-to-noise ratio allows both high probability and a multiplier close to
\(1\).}
\label{fig:single-step-tradeoff}
\end{figure}

\Cref{fig:single-step-empirical} repeats the same calibration across
several signal-to-noise ratios. The empirical exceedance probability tracks the
target probability, and the certified multiplier keeps \(\kappa([B,q])\) close
to \(1\) when the signal-to-noise ratio is high. \Cref{fig:single-step-tradeoff}
plots the certified multiplier
\[
\frac{1+\sqrt{1-\tau^2}}{\tau}
\]
from \Cref{cor-main} against the target probability for several
signal-to-noise ratios.

\begin{figure}[htbp]
\centering
\includegraphics[width=\textwidth]{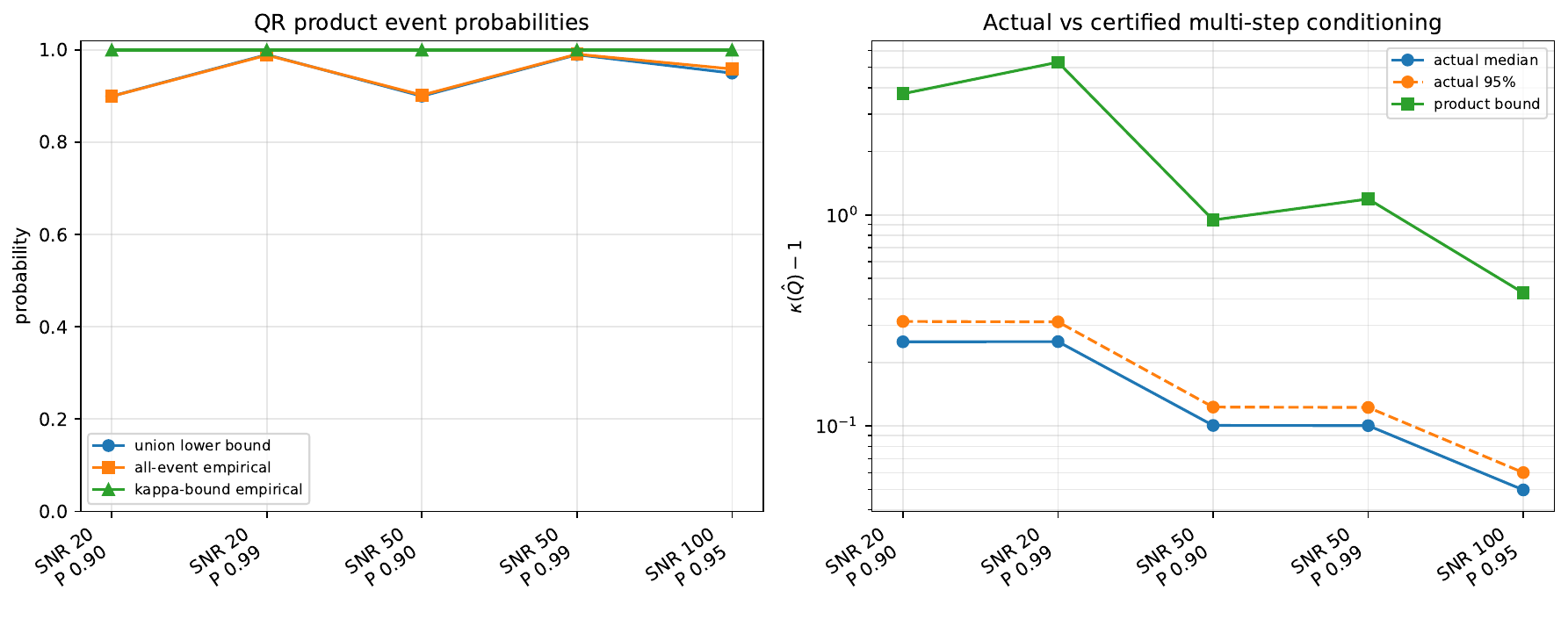}
\caption{Empirical validation of the multi-step product bound
(\Cref{cor:qr-noise-product}). The union bound is conservative, and
the probability that the product condition-number certificate holds is at least
the all-event probability.}
\label{fig:qr-product-empirical}
\end{figure}

For the multi-step process of \Cref{cor:qr-noise-product},
\Cref{fig:qr-product-empirical} compares the union-bound probability and
the product condition-number certificate with simulated QR processes. This
multi-step certificate is unconditional and safe, but conservative: the union
over steps and the product of per-step multipliers accumulate slack. The sharper
calibration is the one-step probability--conditioning tradeoff in
\Cref{fig:single-step-empirical} and~\Cref{fig:single-step-tradeoff}.

\section{Numerical Summary}
\begin{table}[ht]
\centering
\caption{Simulation checks for the stated results. Each row was generated by \texttt{scripts/run\_validation\_simulations.py}.}
\label{tab:validation-summary}
\small
\begin{tabular}{@{}lccc@{}}
\toprule
Claim & Trials / grids & Violations & Diagnostic \\ 
\midrule
Rank-two eigenvalue lemma & 3,000 & 0 & max eig. error 6.22e-15 \\ 
Append-column theorem & 10,597 valid & 0 & median actual/bound 0.985204 \\ 
Min-max factor proposition & 10,000 & 0 & all tested variants passed \\ 
Norm under Gaussian noise & 2 grids & -- & max abs. prob. error 3.549e-03 \\ 
Projected Gaussian noise & 1 grid & -- & max abs. prob. error 2.658e-03 \\ 
Least-squares ratio & 1 grid & -- & max abs. prob. error 1.555e-03 \\ 
QR product implication & 3,000 & 0 & emp. all-events 1.000, union lower 1.000 \\ 
\bottomrule
\end{tabular}
\end{table}

All numerical experiments in this section are reproducible from the Python
scripts available at
\url{https://github.com/alilotfi90/probabilistic-qr-noise}. The simulations are
not used as proofs; they are included to catch algebraic errors, tail-direction
mistakes, and notation mistakes in the probability formulas. Together, the
figures and the table illustrate the main probability--conditioning
tradeoff: higher signal-to-noise ratios allow high-probability guarantees while
keeping the certified condition-number multiplier close to \(1\).

\section{Conclusion and future research}
{\color{black}

We derived deterministic and probabilistic bounds for least-squares residuals, orthogonal projections, and QR-type algorithms under Gaussian perturbations. The deterministic append-column result controls \(\kappa([B,q])\) through the component of the perturbation visible to \(B\). The probabilistic results show that the relevant projected residuals are governed by noncentral chi-square and noncentral-F distributions. These formulas lead to QR condition-number bounds under a model in which normalization is exact and orthogonalization is subject to Gaussian noise.

The primary focus of the paper is numerical linear algebra. Nevertheless, probabilistic perturbation models of this kind may also be useful in settings where orthogonalization or dimension reduction is used inside larger algorithms, such as machine learning and signal processing. QR factorization, singular value decomposition, and principal component analysis are common components of high-dimensional data methods; probabilistic PCA is one example where Gaussian models are already central \cite{tipping1999probabilistic}. Questions of robustness to noise and limited precision also arise in learning systems and large-scale stochastic computation \cite{goodfellow2014explaining,smith2020generalization,collet2008stochastic}. We view these as possible future directions rather than as applications proved in the present paper.

Future work could extend the analysis to non-Gaussian perturbations, to more realistic floating-point error models, or to other algorithms in which an ideal update is orthogonal to a previously computed subspace but the computed update is contaminated by noise.
}

\bibliographystyle{plain}
\bibliography{references}
\end{document}